\documentclass[oneside,english]{amsart}
\usepackage[T1]{fontenc}
\usepackage{epsfig,graphicx,times,hyperref}
\usepackage{babel}
\usepackage{amssymb}

\makeatletter

\providecommand{\tabularnewline}{\\}

 \theoremstyle{plain}    
 \newtheorem{thm}{Theorem}[section]
 \numberwithin{equation}{section} 
 \numberwithin{figure}{section} 
 \theoremstyle{plain}

\makeatother
\begin{document}

\title{A note on Erd\H{o}s-Diophantine graphs and Diophantine carpets}

\author{Axel Kohnert and Sascha Kurz }

\address{University of Bayreuth, Department of Mathematics, D-95440 Bayreuth,
Germany}

\email{\{axel.kohnert\},\{sascha.kurz\}@uni-bayreuth.de}

\begin{abstract}
We give an effective construction for Erd\H{o}s-Diophantine graphs and
characterize the chromatic number of Diophantine carpets. 
\end{abstract}



\maketitle

\section{Introduction}

A Diophantine figure, see i.e. \cite{erdoes_diophantine,gauss_integers,diophantine_figures},
is a set of points on the integer grid $\mathbb{Z}^{2}$ where all
mutual Euclidean distances are integers. We also speak of \emph{Diophantine
graphs}. The vertices are points in $\mathbb{Z}^{2}$ (the coordinates)
and the edges are labeled with the distance between the two adjacent
vertices, which is integral. In this language a Diophantine figure
is a complete Diophantine graph. Two Diophantine graphs are equivalent
if they only differ by translation or rotation of vertices. Due to
a famous theorem of Erd\H{o}s and Anning \cite{ErdoesAnning1} there
are complete Diophantine graphs which are not contained in larger
ones. We call them \emph{Erd\H{o}s-Diophantine graphs}. We will give a
proof of this theorem as we need it for an algorithm later on.

\begin{thm}
\textbf{(Erd\H{o}s, Anning 1945 \cite{ErdoesAnning1})} \label{theorem}

Infinitely many points in the plane with pairwise integral distances
are collinear. 
\end{thm}
\begin{proof}
Let $A$, $B$, and $C$ be three non collinear points and let $k=\max(\overline{AC},\overline{BC})$.
Then there are at most $4(k+1)^{2}$ points $P$ such that the differences
of Euclidean distances $\overline{PA}-\overline{PB}$ and $\overline{PB}-\overline{PC}$
are integers. We see this as follows: Due to the triangle inequality
we have $|\overline{PA}-\overline{PB}|\leq\overline{AB}\leq k$. Thus
$|\overline{PA}-\overline{PB}|\in\{0,1,\dots,k\}$. So $P$ is on
one of $k+1$ hyperbolas. Analog we have that $P$ is situated also
on one of $k+1$ hyperbolas through $B$ and $C$. Because two distinct
hyperbolas intersect in at most $4$ points, there are at most $4(k+1)^{2}$
points $P$. 
\end{proof}
A special class of Diophantine graphs are \emph{Diophantine carpets}
\cite{carpets2,carpets}. These are planar triangulations of a subset
of the integer grid.

\section{Problems}

The authors of \cite{gauss_integers} have posed some open problems
for Diophantine graphs and Diophantine carpets which we would like
to solve in this section.

\subsection{Pythagorean triangles}

Let us denote by $\chi(l)$ the number of all Pythagorean triangles
with hypothenuse $l\in\mathbb{N}$. The question in \cite{gauss_integers}
was to determine the asymptotic of the function $\chi(l)$ when $l\rightarrow\infty$.\\
 Due to Jacobi (1828) we have $\chi(l)=d_{1,4}(l)-d_{3,4}(l)$ where
$d_{r,n}(l)$ denotes the number of divisors (including $1$ and $l$)
of $n$ which are congruent to $r$ modulo $n$, see i.e. \cite{0728.11001}.
So $\chi(l)\in O(n^{\varepsilon})$ for $\varepsilon>0$, see \cite{number_theory_hw}
for a deeper analysis of the divisor function.

\subsection{Erd\H{o}s-Diophantine triangles}

Are there Erd\H{o}s-Diophantine triangles and is there an effective
algorithm to determine all integer points $P$ which extend a given
Diophantine triangle (= complete Diophantine graph of $3$ points)
to a complete Diophantine graph of $4$ points?\\
For such an effective algorithm we can use Theorem \ref{theorem}.
For the given integral points $A=(a_{1},a_{2})$, $B=(b_{1},b_{2})$,
and $C=(c_{1},c_{2})$ we have the following system of equations for
a forth point $P=(x,y)$ being the intersection of two hyperbolas.
\begin{eqnarray*}
\sqrt{(x-a_{1})^{2}+(y-a_{2})^{2}}-\sqrt{(x-c_{1})^{2}+(y-c_{2})^{2}}=d_{1}\\
\sqrt{(x-b_{1})^{2}+(y-b_{2})^{2}}-\sqrt{(x-c_{1})^{2}+(y-c_{2})^{2}}=d_{2}\end{eqnarray*}
 Due to the proof of theorem \ref{theorem} we have $-\overline{AC}\leq d_{1}\leq\overline{AC}$
and $-\overline{BC}\leq d_{2}\leq\overline{BC}$ for $d_{1},d_{2}\in\mathbb{Z}$.
Thus we can solve the corresponding $(2\overline{AC}+1)(2\overline{BC}+1)$
equation systems to determine the possible points $P$.\\
To answer the first question we loop over all Heronian triples, which
are not Pythagorean. (Pythagorean triples can be extended) Heronian
triples are triples of edge lengths, which correspond to an triangle
of rational area. The restriction in search to the Heronian triples
comes from the fact that triangles in the $\mathbb{Z}^{2}$ lattice
are always of rational area as the area is the half of the determinant:\[
\left|\begin{array}{ccc}
a_{1} & a_{2} & 1\\
b_{1} & b_{2} & 1\\
c_{1} & c_{2} & 1\end{array}\right|.\]

Compute in the next step all possible embeddings of such a triple
into $\mathbb{Z}^{2}$ and let an implementation of the above described
algorithm search for possible fourth nodes. If we fail to find a fourth
node we found a Erd\H{o}s-Diophantine triangle. We experimentally
noticed that there are very rare, but we found seven examples with
edge lengths:\[
\begin{array}{c}
(2066,1803,505)\\
(2549,2307,1492)\\
\begin{array}{c}
\begin{array}{c}
(3796,2787,2165)\\
(4083,2425,1706)\\
(4426,2807,1745)\end{array}\\
(4801,2593,2210)\end{array}\\
(4920,4177,985)\end{array}.\]
 This is a complete list of Erd\H{o}s-Diophantine triangles having
an edge of length $\le5000.$\\

\subsection{Further Erd\H{o}s-Diophantine Graphs}

In \cite{erdoes_diophantine} the following two Diophantine figures
were depicted, which the author believed to be Erd\H{o}s-Diophantine
graphs. 

\begin{center}\begin{tabular}{|c|c|}
\hline 
\includegraphics[width=4.4cm]{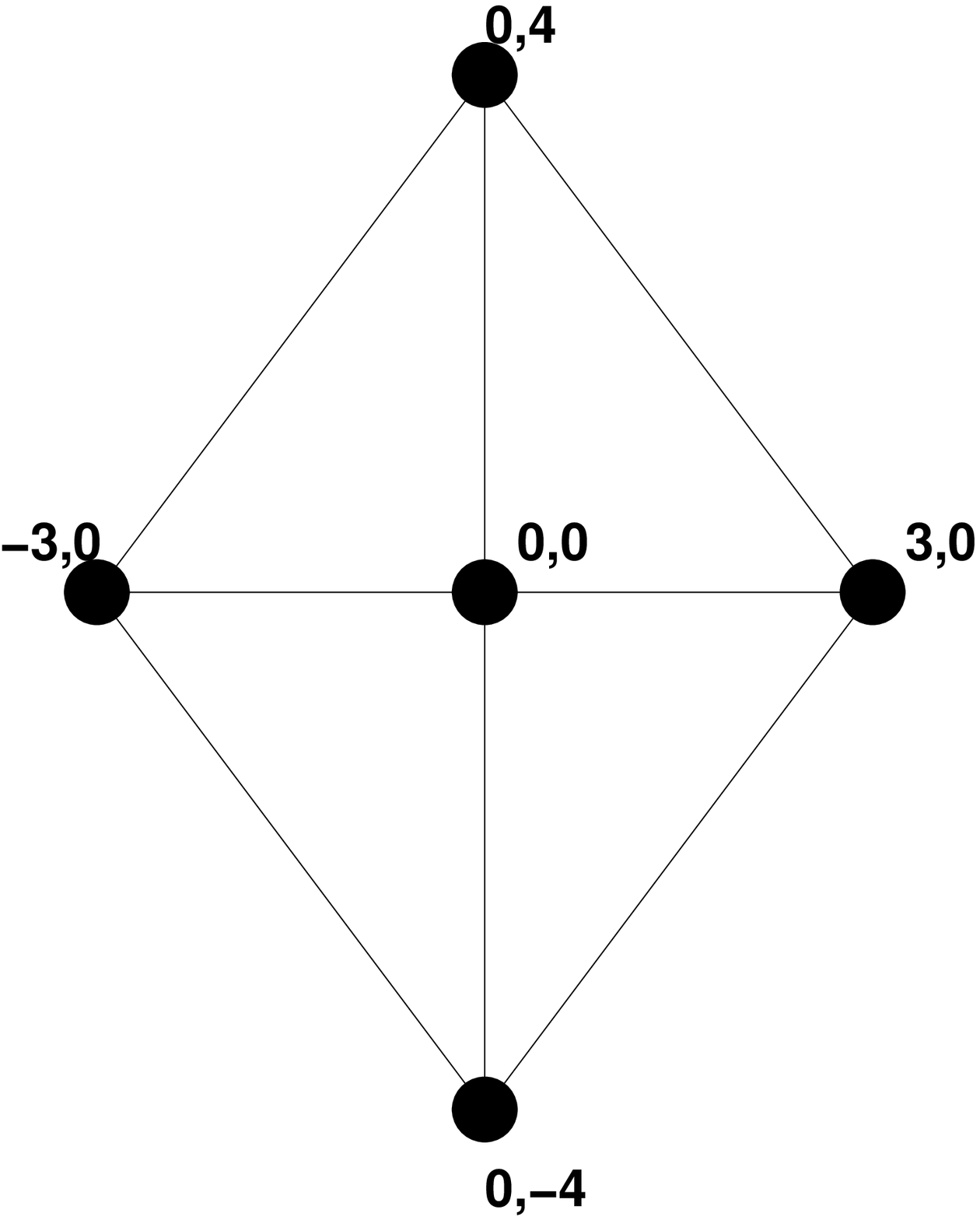} & 
\includegraphics[width=6.0cm]{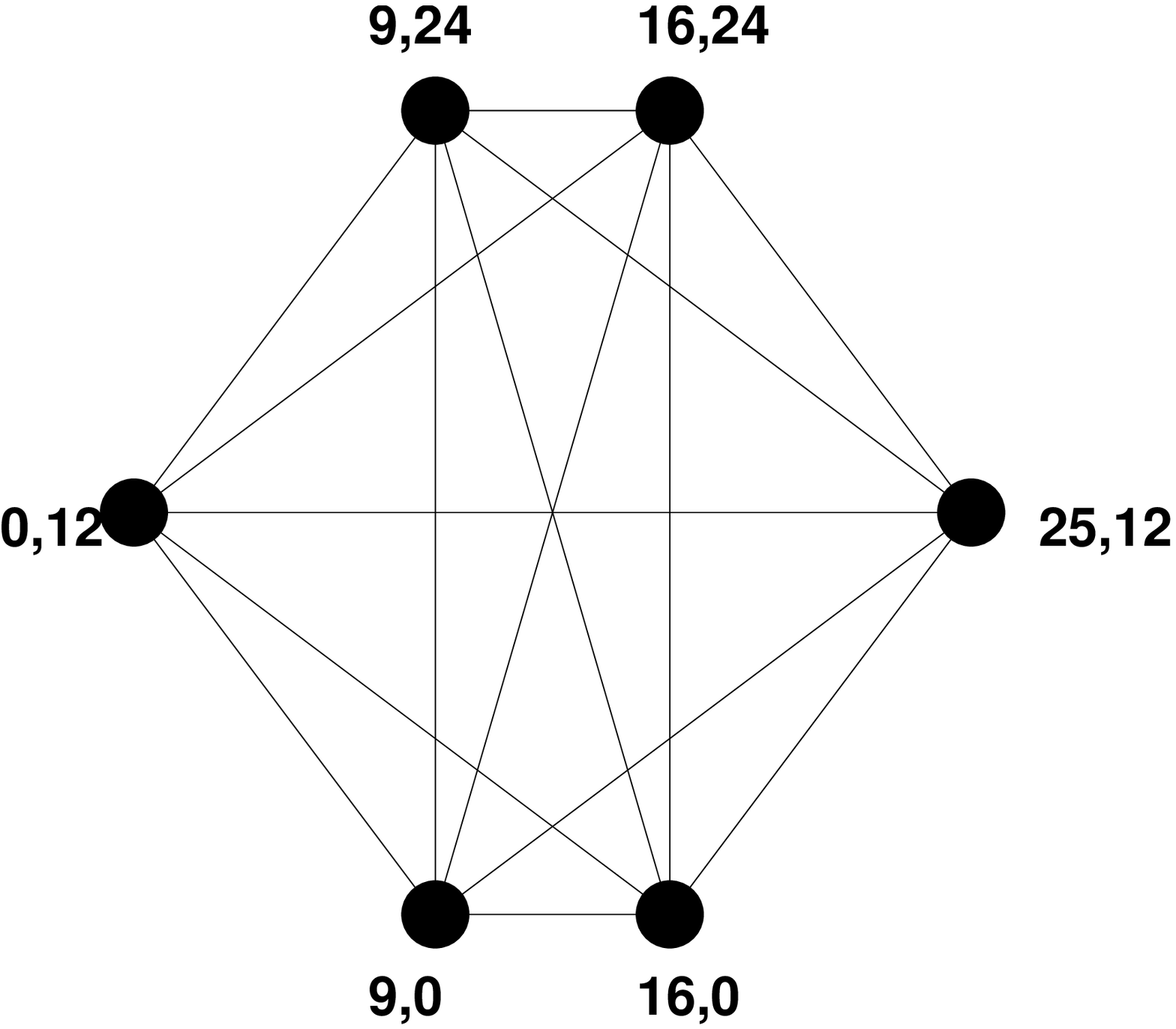}
\tabularnewline
\hline
\end{tabular}\end{center}

With the above algorithm we checked their conjectures, and proved
them.

\subsection{Erd\H{o}s-Diophantine Tetrahedrons}

In \cite{gauss_integers} Pythagorean-Diophantine pyramids were defined
as sets of four points with integer coordinates, integral edge lengths
and three faces being Pythagorean triangles. They asked for Erd\H{o}s
Pythagorean-Diophantine pyramids. We slightly generalize their definition
and search for tetrahedrons with coordinates in $\mathbb{Z}^{3}$
and integral edge lengths, integral face areas and integral volume.
These objects are called \emph{Diophantine tetrahedrons}. In the case
that there is no further point in $\mathbb{Z}^{3}$ having integral
distance to the vertices of the tetrahedron, we call it 
\emph{Erd\H{o}s-Diophantine tetrahedron}.

For four points $A=(a_{1},a_{2},a_{3})$, $B=(b_{1},b_{2},b_{3})$,
$C=(c_{1},c_{2},c_{3})$, and $D=(d_{1},d_{2},d_{3})$ of a Diophantine
tetrahedron we have the following system of equations for a fifth
point $P=(x,y,z)$ being the intersection of three hyperboloids. \begin{eqnarray*}
\sqrt{(x-a_{1})^{2}+(y-a_{2})^{2}+(z-a_{3})^{2}}-\sqrt{(x-d_{1})^{2}+(y-d_{2})^{2}+(z-d_{3})^{2}}=e_{1}\\
\sqrt{(x-b_{1})^{2}+(y-b_{2})^{2}+(z-b_{3})^{2}}-\sqrt{(x-d_{1})^{2}+(y-d_{2})^{2}+(z-d_{3})^{2}}=e_{2}\\
\sqrt{(x-c_{1})^{2}+(y-c_{2})^{2}+(z-c_{3})^{2}}-\sqrt{(x-d_{1})^{2}+(y-d_{2})^{2}+(z-d_{3})^{2}}=e_{3}\end{eqnarray*}
 Using a variant of the above algorithm we did an extensive search
and found several solutions. We give the coordinates of $B,C,D$ where
the first point $A$ is always the origin. We found the following
Erd\H{o}s-Diophantine Tetrahedrons: \\

$\left(\begin{array}{ccc}       


 396 & 132 & 99 \\      

288 & -84 &  0 \\       176 &   0 &  0 \end{array} \right)$, $\left(\begin{array}{ccc}        

             432 & 144 & 108 \\      336 & -48 &  20 \\      297 &   0 &   0 \end{array} \right)$, $\left(\begin{array}{ccc} 

       528 & 396 &  121 \\     468 & 204 & -423 \\     144 & 108 & -135 \end{array} \right),$\\ 

$\left(\begin{array}{ccc}       

             540 & 180 &  135 \\     336 & 252 &    0 \\     400 &   0 &    0 \end{array} \right)$, $\left(\begin{array}{ccc}        

          624 & 468 &    0 \\     648 & 360 & -189 \\     660 & 264 &  -77 \end{array} \right)$, $\left(\begin{array}{ccc}        

            672 & 104 &   0 \\      672 &   0 &   0 \\      600 &   0 & 135 \end{array} \right),$\\

$\left(\begin{array}{ccc} 

          672 &  104 &   0 \\     672 & -104 &   0 \\     600 &    0 & 135 \end{array} \right)$, $\left(\begin{array}{ccc}        

            672 &  153 & 104 \\     672 &    0 & 104 \\     672 &    0 & 0 \end{array} \right)$, $\left(\begin{array}{ccc}  

          672 &  153 & 104 \\     672 & -153 & 104 \\     672 &    0 &   0 \end{array} \right)$.

\subsection{Chromatic number of Diophantine carpets}

We now examine the coloring problem for Diophantine carpets. Clearly
the chromatic number is $1$ iff the carpet consists of union of non-connected
triangles. Given a Diophantine carpet $\mathcal{C}$ we define a graph
$\mathcal{C^{\ast}}$ by replacing the triangles by nodes which are
adjacent iff the corresponding triangles share a common side. This
is the dual graph of the planar graph without a node for the outer
face. As $C$ is a triangulation, each node in $C^{\ast}$ has maximal
degree $3.$ A graph $G$ is bipartite iff it contains no odd cycle
\cite{0945.05002}. For the remaining cases the chromatic number is
$3$. As the chromatic number is according to the theorem of Brooks\cite{0945.05002}
bound by the maximal degree.

\bibliographystyle{abbrv}
\bibliography{erdoes_diophantine_graphs}

\end{document}